\documentclass{ifacconf}

\usepackage{graphicx}      
\usepackage{natbib}        
\usepackage{amssymb}
\usepackage{csquotes}
\usepackage{amsmath}
\usepackage{algorithm}
\usepackage{algpseudocode}
\usepackage{comment}

\newtheorem{assume}{Assumption}
\newtheorem{remark}{Remark}

\newcommand{\algrulehor}[1][.2pt]{\par\vskip.5\baselineskip\hrule height #1\par\vskip.5\baselineskip}
\begin{document}
\begin{frontmatter}

\title{Multi-agent Black-box Optimization using a Bayesian Approach to Alternating Direction Method of Multipliers\thanksref{footnoteinfo}} 

\thanks[footnoteinfo]{D.K. acknowledges support from the  Research Council of Norway through the
IKTPLUSS program and Eindhoven University of Technology. J.A.P. acknowledges support from the National Science Foundation Grant 2036982.}

\author{Dinesh Krishnamoorthy$ ^1 $} 
\author{Joel A. Paulson$ ^2 $}

\address{$ ^1 $Department of Mechanical Engineering,  Eindhoven University of Technology, Eindhoven, 5600MB, The Netherlands. \\
	$ ^2 $Department of Chemical and Biological  Engineering, The Ohio State University,  Columbus, OH 43210, USA.\\
	(email: d.krishnamoorthy@tue.nl, paulson.82@osu.edu)
}

\begin{abstract}                
Bayesian optimization (BO) is a powerful black-box optimization framework that looks to efficiently learn the global optimum of an unknown system by systematically trading-off between exploration and exploitation. However, the use of BO as a tool for coordinated decision-making in multi-agent systems with unknown structure has not been widely studied.
This paper investigates a black-box optimization problem over a multi-agent network coupled via shared variables or constraints,  where each subproblem is formulated as a BO that  uses only its local data. The proposed multi-agent BO (MABO) framework adds a penalty term to traditional BO acquisition functions to account for coupling between the subsystems without data sharing. We derive a suitable form for this penalty term using alternating directions method of multipliers (ADMM), which enables the local decision-making problems to be solved in parallel (and potentially asynchronously).
The effectiveness of the proposed MABO method is demonstrated on an intelligent transport system for fuel efficient vehicle platooning. 
\end{abstract}


\end{frontmatter}

\section{Introduction}

Multi-agent decision-making problems (often formulated as optimization tasks) arise in a wide range of application areas, such as in collaborative robotics, intelligent transport systems, integrated energy systems, sensor networks, and modular manufacturing, where a collection of agents work together to achieve one or more well-defined goals.
A common framework for multi-agent decision-making is the decomposition-coordination approach, where the agents solve their own local subproblems, and a central coordinator  accounts for the coupling/interactions between the different agents, thus facilitating a coordinated action \citep{bertsekas2016nonlinear,lasdon2002optimization}. One such popular distributed optimization approach is the alternating directions method of multipliers (ADMM) \cite{boyd2011ADMM}.  
 
In the traditional optimization framework, each local agent requires analytical expressions (i.e., models) that relate the local decision variables to the local objectives and  constraints in every subsystem. Obtaining reliable models is a time consuming and tedious work, and modeling effort is often one of the biggest bottlenecks in many application areas \citep{chachuat2009adaptation}.
Accurate models may also be challenging to obtain in some applications, due to lack of domain knowledge. This requirement hinders the use of traditional distributed optimization methods that rely on exploiting the structure of such models (in the form of gradients or relaxations) on applications for which accurate models are not readily available. 
 
A class of optimization methods known as ``black box'' optimization is a promising alternative that circumvents the need for detailed mathematical models. Here, the system is considered to be an unknown “black-box” that can be queried through a series of simulations or experiments. The data observed from such evaluations can then be used to find the optimum.  Bayesian optimization (BO) is one such zeroth-order black-box optimization scheme, where the optimum is learned by sequentially interacting with the system \citep{kushner1964new,shahriari2015taking,jones1998efficient}.  

However, in the context of multi-agent systems, formulating the subproblems using BO does not account for the coupling and the interactions between the different agents. In BO, the next sequence of actions are computed by optimizing a so-called \enquote{acquisition function} induced from a probabilistic surrogate model for the local objective function. 
Several different acquisition functions have been proposed in the BO literature, which fundamentally attempt to trade-off exploration (learning for the future) and exploitation (immediately advancing toward goal). Interested readers are referred to the review papers \citep{shahriari2015taking,frazier2018tutorial} for detailed discussions on the different type of acquisition functions and the relationship between them. 
In the distributed setting, however, the local acquisition function in each subsystem does not account for the coupling to the other subsystems. A truly distributed BO algorithm requires the local acquisition functions to be modified so that, e.g., a central coordinator can influence the local BO decision-making in each subsystem. Despite the growing popularity of BO, decomposition-coordination frameworks for BO has received very limited attention in the literature. 

To address this gap, this paper proposes a multi-agent Bayesian optimization (MABO) framework that is capable of solving  black-box optimization problems defined over a multi-agent network coupled via shared variables or constraints. An important and realistic assumption made about the network is that direct data sharing is not allowed between the local agents due to data privacy concerns or bandwidth limitations. Such constraints would not be satisfied by traditional ``centralized'' BO methods, but is naturally handled by our proposed MABO approach, which augments the local acquisition functions with additional penalty terms expressed directly in terms of limited information received from the central coordinator. MABO can be used in combination with any existing acquisition function from the BO literature and, thus, provides a bridge between the developments in BO and distributed optimization. 


The reminder of the paper is organized as follows. Section~2 introduces the preliminaries of Bayesian optimization. Section~3 describes the proposed MABO framework and corresponding penalty terms that are derived using ADMM. The effectiveness of the MABO framework is then demonstrated on a fuel efficient vehicle platooning problem in Section~4, wherein the goal is to minimize overall fuel consumption of the platoon by selecting an optimal cruising speed. 
Lastly, we conclude the paper and discuss some interesting directions for future work in Section~5. 

\section{Preliminaries - Bayesian Optimization} 
Consider the optimization problem 
\begin{equation} \label{eq:min-f}
	\min_{x\in \mathcal{X}} f(x)
\end{equation}
where $ x \in \mathcal{X} $ is the decision variable chosen from a known compact set $ \mathcal{X} $ and $ f: \mathcal{X} \rightarrow \mathbb{R} $ is some cost black-box function whose mathematical structure is unknown.
Bayesian optimization (BO) consists of the following two steps:
\begin{enumerate}
	\item \textbf{Learning}: A probabilistic surrogate model, typically a Gaussian process (GP), is derived for the unknown cost function $ f(x) $. The prior distribution for the surrogate model is updated every time that a noisy observation of the cost $ y = f(x) + w$, where $w$ is a stochastic noise term, is available by conditioning on the observations \citep{rasmussen2006gaussian}.
	\item \textbf{Decision-making:} The posterior surrogate model is then used to induce a suitable acquisition function  $ \alpha: \mathcal{X} \rightarrow \mathbb{R} $, which provides a measure of how desirable querying a point $ x \in \mathcal{X} $ is expected to be with respect to optimizing the cost. Using this, the next action to evaluate is computed by solving\footnote{Note that we choose to define the acquisition function $\alpha(\cdot)$ in a way that its a minimization problem  (as opposed to maximization problem) without loss of generality. }
	\begin{equation}\label{key}
		x^{k+1} = \arg \min_{x\in \mathcal{X}} \; \alpha(x|\mathcal{D}^k),
	\end{equation}
	where $\mathcal{D}^k := \{ (x^j, y^j) \}_{j=1}^k$ is the past set of observations at time $k$ and $ \alpha(x|\mathcal{D}^k) $ denotes the acquisition function evaluated at $x$ induced by the posterior conditioned on $\mathcal{D}_k$.
\end{enumerate} 
\begin{figure}
	\centering
	\includegraphics[width=0.45\linewidth]{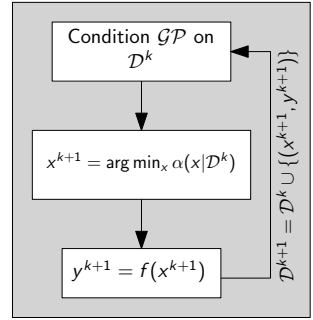}
	\caption{Illustration of Bayesian optimization scheme to minimize an unknown function $ f(x) $.}\label{Fig:BO}
\end{figure}
The choice of acquisition function controls the trade-off between exploration and exploitation when selecting a suitable next action, with the goal being to find the global optimum in as few iterations as possible. 
There are several choices for the acquisition function such as probability of improvement, expected improvement, upper confidence bound, and Thompson sampling. See \cite{shahriari2015taking,frazier2018tutorial} for a more detailed description of the different acquisition functions.

In summary, BO is a black-box decision-making strategy where we sequentially query the system, observe the cost, update our prior belief and choose the next action to query. The BO loop is schematically represented in Fig.~\ref{Fig:BO}.

\section{Proposed Multi-agent Bayesian optimization (MABO) Framework}

Now, let us consider a multi-agent version of \eqref{eq:min-f} that has the following additively separable cost function
\begin{equation}\label{Eq:Additively Separable cost}
	\min_{x\in \mathcal{X}}\; f(x) =  \textstyle\sum_{i=1}^{N}f_{i}(x),
\end{equation}
where $f_i : \mathcal{X} \to \mathbb{R}$ denotes a local cost function for the $i^\text{th}$ subsystem. We make the following assumption about these local cost functions.
\begin{assume} \label{assumption:unknown-local}
	The local cost functions $f_i(x)$ are unknown but can be evaluated at any desired $x \in \mathcal{X}$ for all $i = 1,\ldots, N$. 
\end{assume}
If this was our only assumption, we could apply traditional BO methods (Fig.~\ref{Fig:BO}) on the aggregated set of data from all subsystems. It is worth noting that an improved version of this algorithm can easily be developed by exploiting the additive cost structure in \eqref{Eq:Additively Separable cost}, as discussed in \cite{kudva2022efficient}. However, such methods are centralized in the sense that they require the cost data to be gathered from all local subsystems at every iteration. There are many situations where this is not desired or even possible, which we summarize in the following additional assumption.
\begin{assume} \label{assumption:no-sharing}
    Data cannot be shared across subsystems such that subsystem $i$ cannot directly use data from any other subsystem $j \neq i$. 
\end{assume}


\begin{figure*}
	\centering
	\includegraphics[width=0.8\linewidth]{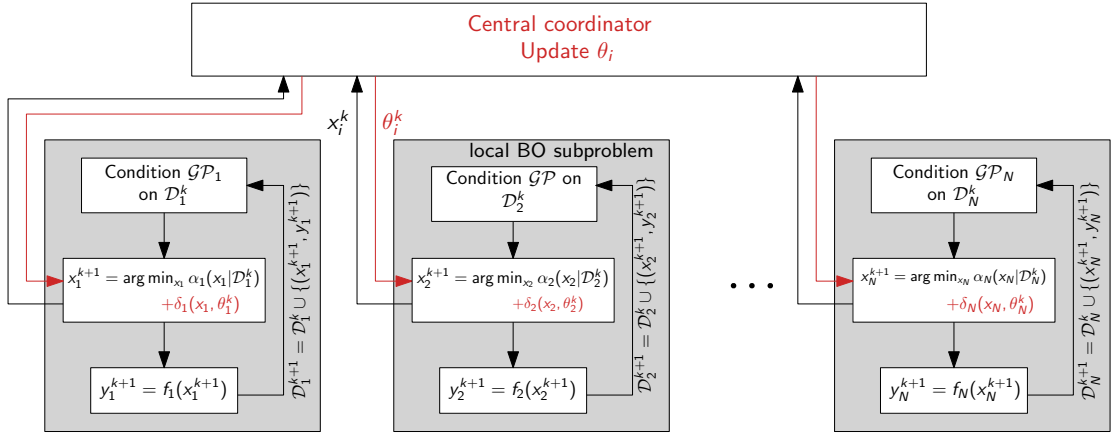}
	\caption{Multi-agent Bayesian optimization (MABO) scheme where each subproblem minimizes its unknown local cost function $ f_{i}(x_{i}) $ using a modified BO method. Each subproblem learns only from the local data set  $\mathcal{D}_{i}$ which is not shared among the different subsystems. The acquisition function in each subproblem has additional penalty terms as a function of the master variables $\theta_{i}$ that enables coordinated action of the different subsystems.  }\label{Fig:FedBO}
\end{figure*}
A straightforward approach to tackle \eqref{Eq:Additively Separable cost}, under Assumptions \ref{assumption:unknown-local}--\ref{assumption:no-sharing}, is to separately apply BO to each subproblem. That is, learn a GP surrogate model for each local cost function $f_i(x_i)$ based only on its local observations $\mathcal{D}_i^k = \{ (x_i^j,y_i^j) \}_{j=1}^k$. Any suitable acquisition function $ \alpha_i(x_{i}|\mathcal{D}^k_{i}) $ induced by the GP posterior is then used for local decision-making by each agent. However, computing the next action by optimizing only the acquisition function induced locally 
 \[ x_{i}^{k+1} = \arg \min_{x_{i}\in \mathcal{X}} \; \alpha_{i}(x_{i}|\mathcal{D}_{i}^k), \]
results in each of the different subsystems taking independent decisions without any coordination. Solving the multi-agent version of the problem \eqref{Eq:Additively Separable cost} requires a coordinated solution in which $x_1^\star = \cdots = x_N^\star$, which is not necessarily achieved with the aforementioned strategy.

To enable the required coordination, we propose to incorporate additional penalty terms to the local acquisition function. 
That is, each subsystem now computes its next action $ x_{i}^{k+1} $ by optimizing
\begin{equation}\label{Eq:DistAcq}
	x_{i}^{k+1} = \arg \min_{x_{i}\in \mathcal{X}} \; \alpha_{i}(x_{i}|\mathcal{D}^k_{i}) + \delta_i(x_{i},\theta_{i}^k),
\end{equation}
where  $ \theta^k_{i}  \in \mathbb{R}^{n_{\theta}}$ denotes the  information that is communicated to the $ i^\text{th} $ subsystem at the $k^\text{th}$ iteration (e.g. by a central coordinator), $ x_{i} \in \mathcal{X}$ denotes the local decision variables,  $ \alpha_{i}(\cdot)$ is any suitable acquisition function used to optimize the unknown local cost $ f_{i}(\cdot) $, and $ \delta_i: \mathcal{X}  \times \mathbb{R}^{n_{\theta}} \rightarrow \mathbb{R} $ is the augmented penalty term that can influence the local decisions. We propose to derive an expression for $\delta_i(\cdot)$ using the alternating direction method of multipliers (ADMM) framework, as described below. 

The optimization problem \eqref{Eq:Additively Separable cost} can be equivalently rewritten by introducing a new variable $ x_{0} $ as 
 \begin{align}\label{Eq:ADMM:ConsensusOpt}
 	\min_{x_{0},\{x_{i}\}\in \mathcal{X}}  \; \textstyle \sum_{i=1}^{N}f_{i}(x_{i}) \;
 	\textup{s.t.} \;  x_{i } =x_{0} \quad \forall i =1,\dots,N.
 \end{align}
The augmented Lagrangian for \eqref{Eq:ADMM:ConsensusOpt} is given by
\begin{equation}\label{Eq:AugmentedLagrangian}
	\min_{ x_{0},\{x_{i}\}\in \mathcal{X}}  \; \sum_{i=1}^{N}\left[f_{i}(x_{i}) +  \lambda_i^{\mathsf{T}}(x_{i} - x_{0}) + \frac{\rho}{2}\|x_{i}-x_{0}\|^2\right],
\end{equation}
where $\lambda = (\lambda_1, \ldots, \lambda_N)$ denotes the set of dual variables (or Lagrange multipliers) associated with the equality consraints in \eqref{Eq:ADMM:ConsensusOpt} and $\rho > 0$ is a penalty parameter. 
The ADMM scheme, which is derived from the augmented Lagrangian in \eqref{Eq:AugmentedLagrangian}, then consists of the following iterations \cite[Chapter 7]{boyd2011ADMM}
\begin{subequations}
\begin{align}
	x_{0}^{k+1} &= \frac{1}{N}\sum_{i=1}^{N}\left[x_{i}^k + \frac{\lambda_{i}^k}{\rho}\right], \label{Eq:ADMM:x0}\\
	x_{i}^{k+1} & = \arg \min_{ x_{i}\in \mathcal{X}} f_{i}(x_{i})+  \lambda_i^{k\mathsf{T}}x_{i} + \frac{\rho}{2}\|x_{i}-x_{0}^{k+1}\|^2, \label{Eq:ADMM:xi}\\
	\lambda_{i}^{k+1} &= \lambda_{i}^k + \rho(x_{i}^{k+1}-x_{0}^{k+1}), \label{Eq:ADMM:lambda}
\end{align}
\end{subequations}
Since $ f_{i}(x_{i}) $ in each subproblem \eqref{Eq:ADMM:xi} is unknown under Assumption \ref{assumption:unknown-local}, we propose to replace it with a local acquisition function $ \alpha_{i}(x_{i}|\mathcal{D}_{i}^k) $. In this case, \eqref{Eq:ADMM:xi} has the same form as \eqref{Eq:DistAcq} with the penalty term equal to
\begin{align}
    \delta_i(x_{i}, \theta_i^k) = \lambda_i^{k\mathsf{T}}x_{i} + \frac{\rho}{2}\|x_{i}-x_{0}^{k+1}\|^2,
\end{align}
where $\theta^k_i := (x^{k+1}_0, \lambda_i^k)$ denotes the coordinating variables, which are updated as shown in \eqref{Eq:ADMM:x0} and \eqref{Eq:ADMM:lambda} at every iteration; these updates are very simple, can be stored with little memory, and ensure that Assumption \ref{assumption:no-sharing} is satisfied. We refer to this overall proposed framework as multi-agent BO (MABO), which is illustrated in Fig.~\ref{Fig:FedBO}. We also provide a complete algorithmic description of MABO in Algorithm \ref{alg:MABO}, which has several unique features that we highlight in the series of remarks below.




 \begin{algorithm}[t]
	\caption{Multi-agent Bayesian optimization.}
	\label{alg:MABO}
	\begin{algorithmic}[1]
		\Require $ N $ local GP models for the local cost functions $ f_{i}(x_{i}) $, initial values $ x_{i}^0 $ and $ \lambda_i^0 $, and penalty $ \rho>0 $
		\algrulehor
		\For {ADMM iterates $k = 0,1,2,\dots$} 
		\State $ x_{0}^{k+1}  \leftarrow \frac{1}{N}\sum_{i=1}^{N}\left[x_{i}^k + \frac{\lambda_{i}^k}{\rho}\right] $ \Comment{Central collector}
		\For {subproblems $ i =1, \dots,N$ } (in parallel)
		\State $ \theta^k_i \leftarrow [x_{0}^{k+1},\lambda_i^k]^{\mathsf{T}} $
		\State $ \delta_i (x_{i},\theta_i^k)\leftarrow \lambda_i^{k\mathsf{T}}x_{i} + \frac{\rho}{2}\|x_{i}-x_{0}^{k+1}\|^2$
		\For {local BO iterates $j = 1,2,\dots$}
		\State Update  local GP model using $ \mathcal{D}_{i}^k $
		\State Induce any acquisition function $ \alpha_{i}(x_{i}|\mathcal{D}_{i}^k)$
		\State $ 	x_{i}^{k+1,j}  \leftarrow \arg \min_{ x_{i}\in \mathcal{X}} \alpha_{i}(x_{i}|\mathcal{D}_{i}^k)+  \delta_i(x_{i},\theta^k_i)$
		\State Query and observe the local cost $ f_{i}(x_{i}^{k+1,j}) $
		\State $ \mathcal{D}_{i}^k  \leftarrow \mathcal{D}_{i}^k  \cup \{(x_{i}^{k+1,j},f_{i}(x_{i}^{k+1,j}))\}$
		\EndFor
		\State $ \lambda_{i}^{k+1} \leftarrow \lambda_{i}^k + \rho(x_{i}^{k+1,j}-x_{0}^{k+1}) $
		\EndFor
		\EndFor
		\algrulehor
		\Ensure Optimal solution $  x_{0} $
	\end{algorithmic}
\end{algorithm}

\begin{remark}[Data privacy]
	The GP model in each subsystem is conditioned only on the local cost measurement, which need not be shared with the other subsystems. This enables a collection of agents to take coordinated actions without sharing its local data with one another, clearly respecting data privacy concerns or any other data sharing limitations (such as bandwidth or connection limitations). 
\end{remark} 

\begin{remark}[Choice of acquisition function]
	Note that the choice of the local acquisition function is independent from the coordinating term $ \delta_i(\cdot) $. Hence any acquisition function proposed in the standard BO framework can be directly used by our MABO approach. Moreover, the choice of the acquisition function can also be different between the different subsystems. Simply put, in a multi-agent network, each agent is free to choose its local acquisition function, without affecting the structure of the MABO strategy.
\end{remark}

\begin{remark}[Multiple BO iterations]
	It can be seen that our proposed decomposable MABO framework (Algorithm \ref{alg:MABO}) has a hierarchical structure with two levels of iterations, namely, the local BO iterations within each subproblem, and the ADMM iterations between the central coordinator and the subproblems. 
	As such, MABO can also be implemented with several BO iterations per ADMM iteration. For example, after every ADMM iterate that updates $ \theta^k_i $, the local subproblems can run one or several local Bayesian optimization iterates for a given $ \theta^k_i $. Furthermore, the different subproblems can also run different number of BO iterations for each ADMM iterate. This is especially beneficial where the local evaluation time and budgets are different in the different agents. This is explicitly shown in Algorithm \ref{alg:MABO}, where the local BO iterations $j$ can be run for any desired finite number of steps. 
\end{remark}

\begin{remark}[Optimal sharing]
Although we have focused our MABO derivation on the consensus problem in \eqref{Eq:ADMM:ConsensusOpt}, the approach extends to \textit{optimal sharing problems} that have coupling constraints of the form $\sum_i A_ix_i = 0$. The main modification needed is to update the form of the penalty term $\delta_i(\cdot)$ whose structure is shown in \citet[Section 7.3]{boyd2011ADMM}. Many important applications involving the distribution of a shared resource would fall into this class of problems (e.g., allocation of energy from a centralized repository in a chemical process). 
\end{remark}

\begin{remark}
Although we considered only the ADMM framework in this section, \eqref{Eq:DistAcq} also allows for other decomposition strategies such as primal decomposition and Lagrangian decomposition, with $\delta_i(\cdot)$ and $\theta_i$ chosen accordingly. 
\end{remark}

\section{Illustrative example: \\ Fuel Efficient Vehicle Platooning}

We demonstrate the performance of MABO on the problem of fuel-efficient vehicle platooning, which is a multi-agent system consisting of several vehicles that take coordinated action to optimize the overall fuel efficiency. Here, ``platooning'' refers to a collection of vehicles that drive in a single-file fashion at close inter-vehicle distance in highways. Platooning of heavy-duty vehicles (e.g., trucks used for freight transport) has several benefits such as reduced traffic congestion, greater fuel economy due to reduced air drag friction, and fewer traffic collisions \citep{franke1995truck,browand2004fuel}. 

Early work on truck platooning focused on tracking inter-vehicle distances using adaptive cruise control \citep{hao2013stability}; however, there has been an increasing interest in optimizing the speed of the platoon to minimize overall fuel consumption. In fact, a recent experimental study showed that, although the fuel economy improved at all platooning speeds, certain speeds resulted in the best overall fuel economy \citep{lammert2014effect}. However, the effect of cruising speed on fuel consumption is not constant and depends on several factors such as freight load, road and driving conditions,  time since last engine maintenance, etc.,  to name a few.
Therefore, developing accurate models that can be used to find the optimal platoon speed is challenging.

An attractive alternative is to resort to a black-box optimization strategy. Black-box formulations are particularly useful when the platoon consists of vehicles that have different configurations and/or manufacturers since it is non-obvious how to coordinate different types of vehicles to minimize overall fuel consumption when each isolated vehicle has a unique optimum speed. 
Such problems can be formulated in the form of \eqref{Eq:Additively Separable cost} where $x_i$ and $f_i(x_i)$, respectively, denote the speed and fuel consumption of the $i^\text{th}$ vehicle. Here, we simulate the fuel consumption of a platoon of $N=7$ vehicles using the approach derived by \cite{sobrino2016reduced}, which leads to a model of the form
\begin{align} \label{eq:fuel-consump}
    f_i(x_i) = a_i + \frac{b_i}{x_i} + c_i x_i + d_i x_i^2,
\end{align}
where speed $x_i$ has units of [km/h] and fuel consumption $f_i(x_i)$ has units of [g/veh-km]. The nominal model parameters are given by the experimentally calibrated values in \citep[Section 5]{sobrino2016reduced}. To emulate $N=7$ completely different trucks, we sample perturbed model parameters from a uniform distribution defined by $\pm 20\%$ of these nominal values.
We highlight that \eqref{eq:fuel-consump} is used only as a source of data generation, as the parameter values and structure are unknown in accordance with Assumption \ref{assumption:unknown-local}.



 \begin{figure*}
	\centering
	\includegraphics[width=\linewidth]{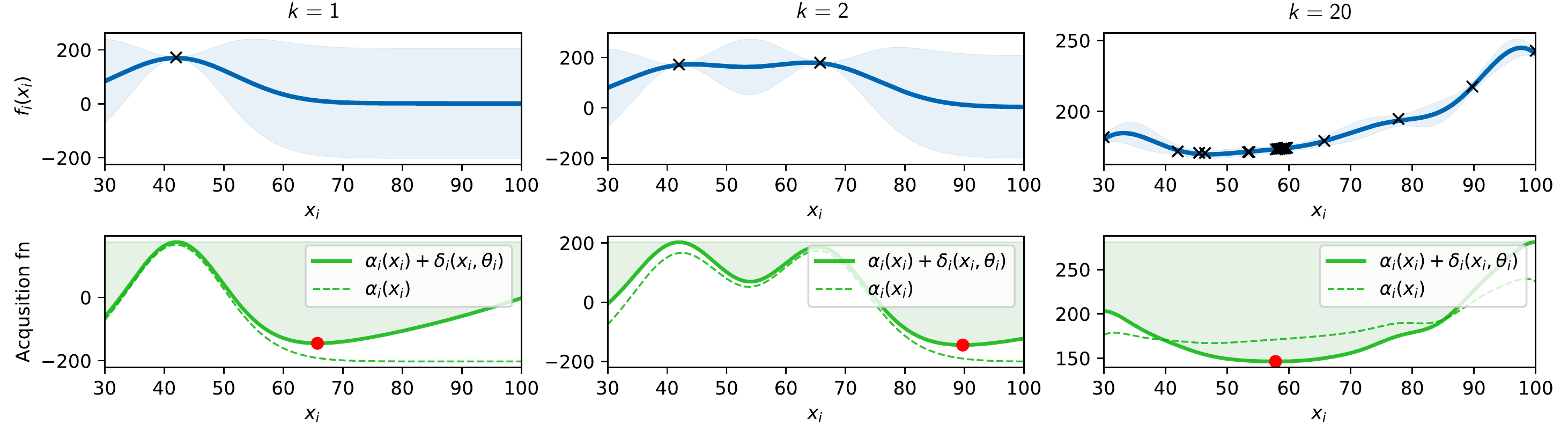}
	\caption{The local Bayesian optimization of one of the subsystems at iterates $ k= 1,2, $ and 20. Top subplots show the Gaussian process of the local cost. Bottom subplots show the  acquisition function with (solid green) and without (dashed green) coordination. The next query point is shown as a red dot.   }\label{Fig:DistBOsteps}
\end{figure*} 


The proposed MABO approach (Algorithm \ref{alg:MABO}) was implemented in \texttt{Python} using the \texttt{GPy} package \citep{gpy2014} to construct and train the GP models using a squared exponential (SE) kernel. We first use a lower confidence bound (LCB) acquisition function defined as follows
\begin{align}
    \alpha_i (x_i | \mathcal{D}_i^k) := \mu_i^k(x_i) - \sqrt{\beta} \sigma_i^k(x_i),
\end{align}
where $\mu_i^k$ and $\sigma_i^k$, respectively, denote the posterior mean and standard deviation predicted by the GP model for the $i^\text{th}$ subsystem at the $k^\text{th}$ iteration and $\beta$ is a scaling parameter used to trade-off between exploration and exploitation. Based on previous work, we select $\beta=4$ throughout this work \citep{srinivas2009gaussian}. We selected a penalty value of $\rho=10$ and only run $j=1$ local BO iteration for every ADMM iteration. 



We first demonstrate the effect of the coordination term in Fig.~\ref{Fig:DistBOsteps}, which plots the GP model and corresponding acquisition functions for a single vehicle at three iterations $k \in \{ 1, 2, 20 \}$. From the bottom subplots, we see that the local non-perturbed acquisition function $ \alpha_{i}(x_{i}|\mathcal{D}_{i}^k) $ (green dashed lines) behaves quite different from the MABO acquisition $ \alpha_{i}(x_{i}|\mathcal{D}_{i}^k)  + \delta_i(x_{i},\theta_i^k) $ (green solid line). Already at iteration $k=1$, the effect of the penalty term $ \delta_i (x_{i},\theta_i^k)$ is evident, which skews the selection process toward lower $x_i^{k+1}$ values such that the next query point is not too far from the average speed of the platoon. By iteration $k=20$, the local BO step converges to around 57 km/hr (near the true optimal value), which significantly differs from the optimum point of the unconstrained acquisition function.


 \begin{figure*}
	\centering
	\includegraphics[width=\linewidth]{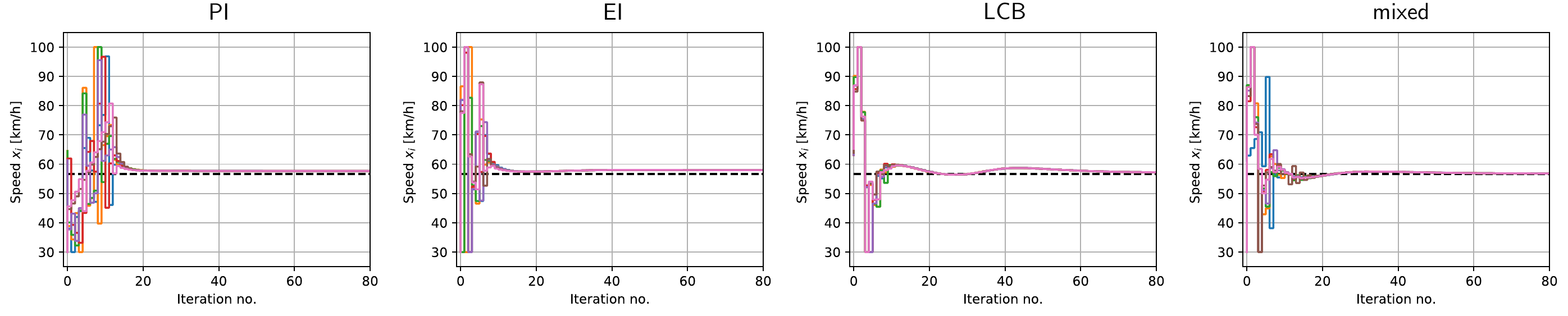}
	\caption{Vehicle speeds when the local BO loops use probability of improvement (leftmost subplot), expected improvement (middle left subplot), LCB (middle right subplot), and mixed (rightmost subplot) acquisition functions. The dashed black line represents the true solution. }\label{Fig:ALLxi}
\end{figure*} 

\begin{figure*}
	\centering
	\includegraphics[width=\linewidth]{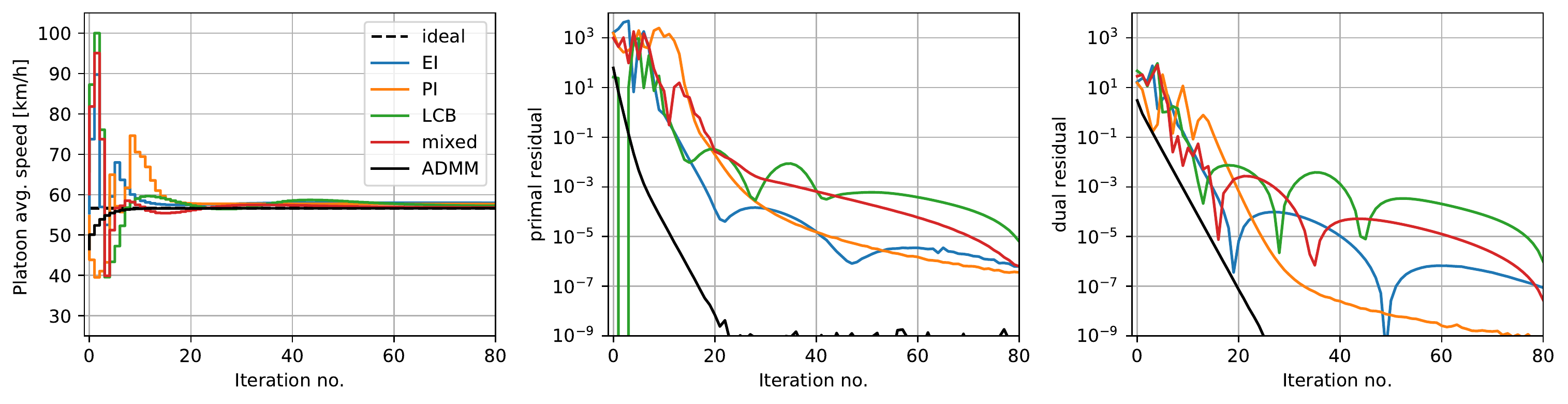}
	\caption{Comparison of platoon average speed, primal residual, and dual residual when local BO loops use probability of improvement (orange), expected improvement (blue), LCB (green), and mixed (red) acquisition functions. The solid black lines represent the ideal ADMM iterations when the local functions are perfectly known. }\label{Fig:ALLx0}
\end{figure*} 

We repeat the same simulation test using two additional types of acquisition functions, probability of improvement (PI) and expected improvement (EI), in addition to LCB. Furthermore, to demonstrate the flexibility of MABO, we also test a ``mixed'' case wherein the different vehicles use different acquisition functions (i.e., vehicles 1 and 2 use EI, vehicles 3 and 4 use PI, vehicles 5 and 6 use LCB, and vehicle 7 using a greedy mean-based policy). The individual vehicle speeds obtained at each MABO iteration using the different acquisition function test cases are shown in Fig.~\ref{Fig:ALLxi}. We also show the true optimum with a black dashed line to serve as a benchmark to compare against.
The average speed of the platoon $ x_{0}^k $, the primal residual $ r^{k} := \sum_i \|x_{i}^k - x_{0}^k\|_{2}^2 $, and the dual residual $ s^k := N\rho^2\|x_{0}^k - x_{0}^{k-1}\|_{2}^2 $  for the different cases are shown in Fig.~\ref{Fig:ALLx0} to compare the effect of the choice of the acquisition function on the progress of the ADMM iterations. It can be seen that the progress of the ADMM iterates are not heavily affected by the choice of the acquisition function. We also benchmark this against the traditional  model-based  ADMM with full model information in each subsystem, which acts as a lower bound on the achievable residuals. This is shown in solid black lines. This empirically quantifies the additional iterations needed to explore the unknown function when using BO in the subsystems. 

\section{Conclusion and Future work}

This paper presents a novel approach for efficiently finding global optimum solutions of multi-agent systems with unknown (black-box) structure and limitations on information sharing by combining elements of distributed optimization with Bayesian optimization (BO). As opposed to independently applying BO in each of the subsystems, the key idea is to modify the local acquisition (or expected utility) function with additional penalty terms that are properly set according to a master coordinator, which requires limited communication bandwidth and memory storage. Our proposed multi-agent BO (MABO) framework is broadly applicable in the sense that it works for different network structures (e.g., consensus and sharing formulations) and allows for the use of any well-defined acquisition function within the local subproblems. Furthermore, the local subproblems can be straightforwardly formulated using constrained BO \citep{gelbart2014constrained,paulson2022cobalt} or safe BO methods \citep{berkenkamp2021bayesian,DK2022SafeBayesOpt} whenever local black-box constraints are present in a given subsystem. The proposed MABO approach is demonstrated on a vehicle platooning problem, where the goal is to select a common platoon speed for multiple (different configuration) vehicles such that the overall fuel consumption is minimized. We show that MABO is capable of quickly identifying the optimal platoon speed, even with no prior knowledge of the relationship between speed and fuel consumption in any of the vehicles in the platoon. This work represents a first step toward a decomposable BO framework and future work should involve better understanding the theoretical properties (e.g., rate of convergence) of the proposed method. In addition, it would be interesting to explore the impact of the update rule chosen for the central coordinator on performance, as there may exist strategies that accelerate the speed of convergence in certain applications. 


\bibliography{L4DC,BibOpt}             
\end{document}